\documentclass[12pt,a4paper]{article}
\usepackage[utf8]{inputenc}
\usepackage{amsmath}
\usepackage{amsfonts}
\usepackage{amssymb}
\usepackage{url}
\usepackage{graphicx}

\usepackage{epstopdf}
\AppendGraphicsExtensions{.pdf}
\graphicspath{ {images/} }

\newtheorem{remark}{Remark}
\newtheorem{lemma}{Lemma}
\newtheorem{corollary}{Corollary}

\newtheorem{algorithm}{Algorithm}

\author{Alexander Kushkuley, kushkuley@gmail.com  \\
Joshua Correa, jcorrea@salesforce.com  }
\title{Improving Recommendation Relevance by simulating User Interest   }
\begin{document}
	\maketitle

\begin{abstract}
\noindent Most if not all on-line item-to-item recommendation systems  rely on estimation of a  
distance like measure (rank) of similarity between items.  For on-line recommendation systems, time sensitivity of this similarity  measure is extremely  important.  We observe that recommendation "recency"  can be straightforwardly and transparently maintained by iterative reduction of ranks of inactive items. The paper briefly summarizes algorithmic developments based on this self-explanatory observation. The basic idea behind this work is patented in a   context of online recommendation systems.   
\end{abstract}

\maketitle

\section{Introduction}
Most if not all on-line item-to-item recommendation systems  (cf. e.g. \cite{Wic1}, \cite{Wic}) are relying on a notion of product similarity so that for an (anchor) item $A$, $A$-recommended items $ B_1, B_2, \cdots B_m $
are displayed on the products $A$ web page in the order of their similarity to $ A $. For further reference denote such a page as 
\begin{equation}
	A \rightarrow  (B_1,r_1), (B_2,r_2) \cdots , (B_m, r_m) 
\end{equation}
where $ \{ r_1, r_2, \cdots , r_n \}  $ are some numbers (ranks) associated with a list of items $ L = \{B_i, \; i = 1, \cdots, m  \} $ .
Usually, the ranks $ r_i $ are obtained via an estimation of some pair-wise similarity measure $\mu$ that  is defined on the set of all product pairs so that by
definition (see e.g. \cite{Amazon})
\begin{equation} 
	r_i  =   \mu(A,B_i) , \;  i = 1, \cdots, m 
\end{equation}
On-line  product recommendations are highly time-sensitive and various techniques are employed 
to keep "distance"  measure $\mu$ time relevant
(see e.g \cite{Amazon}).  For example,  some of these  techniques could be based on event frequency estimation  (cf. e.g. \cite{HH1}, \cite{HH}, \cite{HK})),  some other time aware methods employ random walks (cf. e.g. \cite{Nik}).   In some cases  ad hoc procedures that take event time-stamps into account are implemented.    
\newline \newline We propose a simple  event-log simulation method that 
\begin{itemize}
	\item [-]  creates recommendation lists (1) while avoiding pair-wise similarity measure estimation (2)
	\item   [-] maintains items "recency" by iteratively reducing ranks of inactive items, thus avoiding frequency counting
	\item [-] directly takes existing recommendation clicks into account 
	\item  [-] does not use any time-dependent logic to adjust recommendation ranks 
\end{itemize}

\section{Description of the  Algorithm} 
Let
\begin{equation}
	q_i = \mathbb{P}( B_i | \textrm{recommendation click on page A)} , \;  i = 1, \cdots, m
\end{equation}
be probabilities of recommendation clicks on products $B_i$  under condition that a user is visiting a web page   corresponding to a product $A$ .  
\begin{remark}
	Ranks $r_i $ in (1-2) can be viewed as  an estimate of probability  distribution (3), since   
	the list of ranks  in (1)   can be turned into  a  list of probabilities  $ Q = \{    q_1, q_2 \cdots, q_m \} $, for example,  by setting $  q_i = (1 - q'_i)/(m-1)$ where $ q'_i = r_i / \sum_{j=1}^m r_j  $
\end{remark}
\noindent Of course, the probabilities $q_i $ in (3) are not known (at least in  real time) and standard way of estimating distribution (3) is event frequency counting (cf. e.g.  \cite{HH1}, \cite{HH}). However, plain frequency counting in a context of on-line  recommendation systems has a tendency to overestimate items that were popular in the past  and requires complicated logical rules for counter adjustments
(cf. nevertheless \cite{HK}). Our goal, therefore,  is to design a simple on-line  
click-probability estimation algorithm with  built in time sensitivity. Suppose that we have a categorical distribution 
\begin{equation}
	P = \{p_1, p_2 \cdots, p_m \}
\end{equation}
that corresponds (as in (1) and(3))  to a list of items
\begin{equation}
	A \rightarrow  L =  \{ B_1,B_2, \cdots, B_m \} 
\end{equation} 
Note, that a  click on a recommendation $B_i$ can be represented by a delta function distribution
\begin{equation}
	\delta_i = (0, \cdots , 1 ,   \cdots , 0)  
\end{equation}
on  a set of indexes  $ \{1,\cdots,m\} $ that has all its mass values equal to 0 except the one with index   $i.$ It is thus quite natural  to assume that such a click triggers convex mixture  update of probability distribution $P$  by a rule 
\begin{equation}
	P \rightarrow  \alpha P + ( 1 - \alpha) \delta_i  , \; \; 0 < \alpha < 1  \nonumber
\end{equation}  In other words, a click on a recommendation  $B_i$ reduces ranks of all other page recommendations, while tilting estimated  rank-distribution towards recommendation of the item $B_i$ in a simplest way possible.    
Thus we arrive at the following online recommendation algorithm.
\begin{algorithm}
	
	Fix a close to 1 real number $ 0 < \alpha < 1 $ (called bellow rank reduction parameter)  and  a small threshold $\epsilon$ (e.g $\epsilon=.001$). 
	A click on a recommendation  $B$ with an anchor $A$ causes the following transformation of the categorical distribution (4) associated with the list  (5) :
	\begin{enumerate}
		\item [1)] if 	$ B = B_j \in L $, set
		$ p_i \rightarrow \alpha  p_i, \;  i = 1, \cdots n, \; i \neq j $ and set 
		$ p_j  \rightarrow  \alpha p_j + 1-\alpha $ 
		\item [2)] if $P = \emptyset $ then  set the list of items recommended by $A$  to $\{B\}$  and set $ P = \{1\} $
		\item [3)]  if  $ B \notin L $ then it can be added to the list  with   probability  that is less than $ min \{ p_i, i = 1, \cdots, m \} $  or, what seems to be a better strategy,  a new item is added  by computing a max-entropy mixture distribution as explained below in Appendix 1 (10-11)
		\item [4)]  any item   $ B_i \in L $ such that  $ p_i < \epsilon $  is  removed  from $ L $ and resulting list is rescaled so that  remaining probabilities sum-up to one.  
	\end{enumerate} 
	
\end{algorithm}
Two "built in" features of this algorithm are worth mentioning:
\begin{itemize}
	\item  decreasing $\alpha $ gives higher priority  to recent click recommendation events 
	\item  and therefore, sensitivity of this ranking scheme to new recommendation events can be easily controlled (even at runtime) by adjusting just one parameter 
\end{itemize}

\begin{remark}
	(cf. e.g. \cite{rank}). Suppose that it is desirable that an item should loose half of its rank if it was idle while a list it belongs to  was updated $T$ times.
	That can be achieved by setting parameter $ \alpha $  to   $  \exp( -\log(2)/T) $. For example, if $T = 10 $ then $ \alpha \approx .93 $ 
\end{remark}

\begin{remark}
	Somewhat similar approach  of "counting-with-exponential-decay strategy" in a context of  Heavy Hitters (Elephant Flows) estimation was recently introduced in \cite{HK}.
\end{remark}
\noindent For some intuitive "theoretical"  considerations in favor of Algorithm 1 see \cite{AK} and Appendix 2. A simulation model for evaluation of this algorithm is described in Appendix 3. A few useful properties of the algorithm can be summarized as follows:
\begin{itemize}
	\item [-] An unpopular (idle)  recommendation will loose its rank-probability at an exponential rate 
	$ \alpha^t $ (for any $ t $-clicks missing the item)
	\item  [-] Clicked recommendation probability $ p $
	is scaled up by a multiplier 
	\begin{equation}
		\alpha + p^{-1} ( 1 - \alpha ) > 1    \nonumber
	\end{equation} 
	with magnitude that is inversely related to its  probability.  Therefore a relative  after-click boost for a low rank recommendation is  higher than for a popular one  or, in other words, the algorithm tends to  increase  visibility for less up-till-now  visited items
	\item [-] As was already mentioned, the rates of "idleness decay" and "activity boost" are both controlled by just one parameter
	\item [-] No time dependent logic is executed anywhere
	\item [-] The algorithm can be run in parallel with or as an enhancement of any other recommender(s) deployed on a (web) site 
\end{itemize}

\subsection{Event Log Simulation }
Unfortunately, for large on-line systems a real time implementation of the Algorithm 1 could be complicated due to necessity of simultaneous updates of many items that  change asynchronously. Therefore,  as a practical matter,  a semi-online simulation version of this algorithm is implemented as follows:
\begin{algorithm}

	At regular time intervals (e.g. every hour) 
	the Algorithm 1 is applied
	 to a collection  of lists (5) by simulating all the recommendation click events in the latest unprocessed event log. The resulting collection of probability distributions  is recorded and reused in the next simulation run.
\end{algorithm}

\subsection{Rank Reduction  Parameter Selection}
There are various ways to select parameter $\alpha$. 	One of them, that is currently implemented, uses the following data:
\begin{itemize}
	\item  [-] A time "period of interest" $T$ (e.g. one week) that is selected with an understanding that   a rank of an item that was idle for 
	time $T$  will be significantly reduced (e.g. halved)
	\item [-] current average per-product click rate that is computed from the event log  
\end{itemize}
With this data at hand the parameter  $\alpha$ can be (re)calculated on each simulation run (cf. e.g. \cite{rank}  and Remark 2)

\subsection{Two Enhancements}
Our first observation is that a product could be recommended on multiple web pages, that is, a product $X$ can be recommended by a few products $A,B,C...$. This leads to a following modification of Algorithm 1
\begin{algorithm}
	When  recommendation  $B$ is  clicked, apply Algorithm 1  to all the recorded recommendation lists  $ A' \rightarrow L'  $ such that $ B \in L' $.
\end{algorithm}  
\begin{remark}
	As a modification to the above, one might consider filtering anchors $A'$ by some criteria of closeness to $A$ 
\end{remark}
\noindent It also makes sense  to consider other events besides recommendation clicks.
Important events that are usually tracked by on-line retail systems are, of course, 
check out events  and "add to cart" events. When a reasonable recommender is  running, these events will be less frequent than recommendation clicks. In this case a following modification of Algoritms 1 and 2 could be considered 
\begin{algorithm}.
	\begin{enumerate}
		\item [1)] At simulation start-up, collect and maintain lists of recently checked out items and added to cart items 
		\item [2)]  Compute rank reduction parameters $\alpha_R, \alpha_C, \alpha_A$ for recommendation clicks, check-outs and add to cart events separately (based on their respective click rates)  as described in section 2.2
		\item [3)]  Use appropriate rank reduction parameters computed on step 2) in cases when a recommended item was  recently checked-out or added to cart 
	\end{enumerate}
	
\end{algorithm}  

\section{Conclusion}
Recommender implementations based on  Algorithm 2 and its modifications achieved statistically significant increase of recommendation click rates in a few A/B tests. Relatively recent A/B test showed significant increase of recommendation click rates over baseline for 15 out of 22 participating web sites. It should be noted, however, that Algorithm 2 as a simulation technique depends on at least one other recommender running on the same web site to bootstrap the simulation.

\section{Appendix 1.  Distribution Mixture and Max Entropy Principle}
For any two discrete distributions $ P $  and $ Q $ of the same size, we have a one-parameter  family of convex distribution mixtures 
\begin{equation}
	P' = P^{\prime}(\alpha) = (1 - \alpha )P + \alpha Q, \;\; 0 \leq \alpha \leq 1    \nonumber
\end{equation}
According to Maximum Entropy  principle, (see e.g. \cite{ent}), an unbiased convex mixture of distributions $P$ and $Q$  corresponds to a choice of parameter $\alpha$ that maximizes mixture entropy. As a function of $\alpha$ the entropy can be written down as  
\begin{equation}
	\mathbb{E}(P') = \mathbb{E}(P,Q,\alpha ) = -[\sum_i (1 - \alpha )P_i + \alpha Q_i ]\log\;[\sum_i (1 - \alpha )P_i + \alpha Q_i]   
\end{equation}  
When $Q$ is a delta-distribution (6), the distribution  $P^{\prime} $
differs significantly  from $ P $ in the $i$-th place only, while  its mass at any other index  $ j \neq i $ is just  $ \alpha p_j $ .  Hence, in this case  the  entropy function $ E(P^{\prime}) \equiv E(\alpha) $ looks especially simple
\begin{eqnarray}
	-E(\alpha) = (1 - \alpha q ) \log( 1 - \alpha q ) + 
	\alpha \sum_{j \neq i } p_j \log(p_j) + \alpha\log(\alpha) q   \nonumber  \\
	\text{where } q = 1 - p_i = \sum_{j \neq i } p_j  \nonumber
\end{eqnarray}
It is now clear that  maximization problem $ \max_{\alpha} E(\alpha)  $
has a simple closed form solution. Indeed, the derivative (in $\alpha$)  of $ -E(\alpha)$
is 
\begin{equation}
	-q \log( 1 - \alpha q) - q + \sum_{j \neq i } p_j \log(p_j) + ( 1 + \log(\alpha) ) q  \nonumber
\end{equation}
Setting  this derivative to zero leads to the equation  
\begin{equation}
	q \log \left( \frac{\alpha}{1 - \alpha q } \right) = -C_i   \; \equiv \; -\sum_{j \neq i } p_j \log(p_j) \nonumber
\end{equation}
and we get the following
\begin{lemma}
	In case of delta-function distribution (6), convex mixture distribution (7)
	has maximum entropy when   
	\begin{eqnarray}
		\alpha = \frac{ 1}{ q + e^{C_i/q }}  \nonumber  \\
		p_i^{\prime} = 1 - \alpha q = \frac{e^{C_i/q}}{q + e^{C_i/q }}  \nonumber
	\end{eqnarray} 
	where $ C_i= \sum_{j \neq i } p_j \log(p_j) $  and 
	$ q = 1 - p_i = \sum_{j \neq i } p_j $.
\end{lemma}
\noindent In a case when a clicked item is not in the list  (1) we can pad our distributions $ P $ and $ Q $ like this
\begin{eqnarray}
	P = (0, \; p_1, \; \cdots ,\; p_n)  \\
	Q = (1,\;  \cdots , \;  0 , \; \cdots , 0)  
\end{eqnarray}
and compute the mixture parameter $ \alpha $ in the same way as above, getting
\begin{corollary} An unbiased addition of a so far unseen recommendation click is a convex mixture of distributions (8) and (9)  with a choice of mixture parameter given by
	\begin{eqnarray}
		\alpha = \frac{ 1}{ 1 + e^{H}}   \nonumber \\   
		p_i^{\prime} = 1 - \alpha  = \frac{e^{H}}{1 + e^{H}} \nonumber 
	\end{eqnarray}
	where $ -H $ is entropy of $P$
\end{corollary} 
\begin{remark}
	Computation in Lemma 1 is  primarily motivated by the following considerations.
	On the one hand, there is some empirical evidence confirming that in the case 3) of  Algorithm 1,
	distribution updates suggested by Corollary 1 
	produce better results than
	simple strategy of  placing a new item at the bottom of a list. On the other had, since "correct rank" of a so far unseen  item cannot be compared to the ranks of already visited items, an unbiased max-entropy solution seems to be a good choice. 
\end{remark}	

\section{Appendix 2. Some Theoretical Considerations}
This short section closely follows \cite{AK}. In general, arrival of recommendation clicks  can be viewed  as a stochastic processes (sequence) of discrete time-dependent categorical distributions $ P(\tau) $ where time $\tau$ is measured in number of clicks per web page.     
It is reasonable to assume that distribution $ P(\tau) $ does not change much during small time intervals (of, say a few hundred clicks).   
To evaluate Algorithm 1 ability to  track  recent event frequencies,
consider a following simplified scenario.
Let's assume,  that  number of iterations $  t $  corresponds  to recommendation relevancy time window. For example, if only last week recommendation ranks  are of high importance, let  $t$ be  a
"weekful" of clicks.  Measuring time by click-counter, suppose that estimated click-distribution at the start of the time period was  $X$ and that for time $ t_1 $ the (unknown) incoming click distribution $ P_1 $ did not change. Suppose also that at time $t_1$    the incoming distribution switched to $ P_2  $ and did not change for the remaining time  $ t_2 = t  -   t_1  $.
It is not hard to see (cf. \cite{AK}) that at the end of  time period the expectation of a  distribution computed by Algorithm 1 will be 
\begin{equation}
	\alpha^{t} X + \alpha^{t_2}(1 - \alpha^{t_1}) P_1 + ( 1 - \alpha^{t_2}) P_2  \nonumber
\end{equation} 
For sufficiently large $t$ only the last two terms are significant and the above formula essentially shows that our simulation algorithm approximates  incoming click distribution as a weighted sum of $P_1$ and $P_2$.  To see how this  approximation is affected by recent click events  let's estimate  the ratio of coefficients at $P_2$ and $ P_1 $ in the expression above.
Since $ \beta = 1 - \alpha $ is supposed to be small, we have (cf. \cite{AK})
\begin{equation}
	\frac  {1 - \alpha^{t_2} } {   \alpha^{t_2}(   1 - \alpha^{t_1}   )  }     = 
	\alpha^{-t_2}   \frac   { 1 - (1 - \beta)^{t_2}   }  { 1 - ( 1 - \beta )^{t_1}  } \approx \alpha^{-t_2}  \frac{t_2}{t_1} 
\end{equation} 
\begin{remark}
As we see,	Algorithm 1 introduces approximately  $\alpha^{-1}$ per iteration ”velocity boost” for recently active items 
\end{remark}
\begin{remark}
	It should be clear (cf. \cite{AK}), that in case of $m$  incoming click distributions $ P_1, P_2, \cdots, P_m $ that are "switched over" at times $0 < t_1 < t_2 < \cdots < t_m < t $ the expectation of the distribution produced by Algorithm 1 will be 
	\begin{equation}
	 E(\alpha)= \alpha^{t} X + \alpha^{t-t_1} (1 - \alpha^{t_1}) P_1 +
	 \alpha^{t-t_1 - t_2} ( 1 - \alpha^{t_2}) P_2 + \cdots +
	 ( 1 - \alpha^{t_m})P_m
	   \nonumber 	
	\end{equation}
Dividing by the coefficient at $ P_m$, throwing away the 
first term  (as $ \alpha^{t} \ll 1 - \alpha^{t_m} $) and applying the same approximation as in (10), one gets 
\begin{equation}
	 E(\alpha) \approx \;\alpha^{t-t_1} \frac{t_1}{t_m} P_1 \; + \;   \alpha^{t-t_1 - t_2}  \frac{t_2}{t_m}P_2 \; +  \; \cdots \; + \; P_m   \nonumber
\end{equation}
The last term dominates this expression if ratios $ t_i/t_m$ are not too large for all $ i = 1,2,\cdots,m$. In other words, under  described conditions, timeliness  of the Algorithm 1 approximation largely depends on the incoming distribution switching rate.

\end{remark}
\section{Appendix 3. Simulation of a Random Walk on a Standard Simplex }
 In general, a problem of estimating  approximations produced by the algorithms described in Section 2  seems to be hard for general arrival processes. Hence, we describe here a realistic  simulation model that is based on the preceding discussion. 
 Let's  consider  the space of all discrete $n$-dimensional distributions, i. e. the standard $n$-dimensional  simplex $\Sigma^{n-1} \subset \mathbb{R}^n $. A random walk that we would like to simulate starts at the center $ X_1 = (1/n, 1/n, \cdots, 1/n) $ of $\Sigma^{n-1}$. This choice corresponds to an assumption that a recommended  item initially has  $n$-recommendations and that initially all of them are indistinguishable, having the same rank,  
 $R_1 = (1,1,\cdots, 1)$.  An incoming click is thus represented by a selection of a vertex of $\Sigma^{n-1}$ as described in Section 2 (cf. (6)).  
 \subsection{Incoming Click Distributions
 }The first moving part of our simulation model is a  sequence 
 $C_1, \cdots, C_{m-1}  $ of categorical 
  distributions represented by a sequence of points in $\Sigma^n$.  
 It is  assumed that the total number of simulated clicks $T$ is split into intervals (days) of random length $ t_1, t_2, \cdots , t_m, \;\; \sum_{i}^mt_i = T $. Based on the considerations discussed in  Remark 7, let us assume that  click distribution changes occur at times
 $t_1, \; t_1 + t_2, \; \cdots, \; t_1 + t_2 \cdots + t_{m-1}$. In other words, 
 for any $ t $ in the semi-open interval $ T_k \equiv  ( t_1 + t_2 + \cdots + t_k, \; t_1 + t_2 + \cdots + t_{k+1}  ] $, the probability of an item $i$ to be clicked at "time" $t$ (that is for the i-th vertex of $\Sigma^{n-1}$ to be selected) is equal to the $i$-th coordinate of $C_k$. For simplicity, and again in  accordance with the Remark 7, we will assume that interval lengths $ |T_k|$ are sampled from a normal distributed $N(\mu, \sigma)$ with $ \sigma=\frac{1}{5}\mu $.  Further,  
  it is reasonable to assume that the point $C_{i+1}$ is close the point $C_i$, for example, $C_{i+1}$  differs from $C_i$ by a small Gaussian noise. The following Julia code (cf. \cite{Distr}) illustrates click distribution switching procedure that satisfies described assumptions.
 \begin{verbatim}
 using Random, Distributions, LinearAlgebra	
 function nextDistribution(C::Distributions.Categorical; eps=.01 )
 	 X = C.p                        # previous categorical distribution
 	 n = size(X)[1]                 # dimension of the simplex
 	 S = rand(n,n); S = S'S         # random symmetric positive definite matrix
 	 I_n = Diagonal( ones(n) )
 	 S = I_n + eps * S              # small positive definite perturbation of I_n 
 	 G = MvNormal(X, S)             # Gaussian with the mean X and covariance S
 	 V = abs.( rand(G, 1) )[:,1]    # "almost G-distributed" positive vector
 	 X = V/sum(V)                   # next categorical distribution
 	 return Categorical(X)
 end
 \end{verbatim}
\subsection{Simulation of Algorithm 1}
The second moving part of the simulation is  a straightforward implementation of Algorithm 1  applied to just one anchor, with the assumption that item clicks come from the sequence of categorical distributions $C_k, \; k = 1,2, \cdots $
Let  $ t+1 \in T_k $, let $ i $ be a random index sampled from the distribution $C_k$ at time $t +1 $ and let $\alpha$ be a fixed smoothing parameter from the description of Algorithm 1 (cf. Section 2.2). Then  
\begin{equation}
	X_{t+1,j} = \alpha X_{t,j} \;\; \text{for} \;\; j= 1,2, \cdots n, \; j \neq i  \nonumber 
\end{equation}
and
\begin{eqnarray}
	X_{t+1,i} = \alpha X_{t,i} + 1 - \alpha  \nonumber \\
	R_{t+1, i } =  R_{t+1, i } + 1  	\nonumber
\end{eqnarray}
where for all variables, the first index is the simulation iteration number (click "time") and the second index is a vector coordinate.

\begin{figure}
	\centering
	\includegraphics[width=1\linewidth]{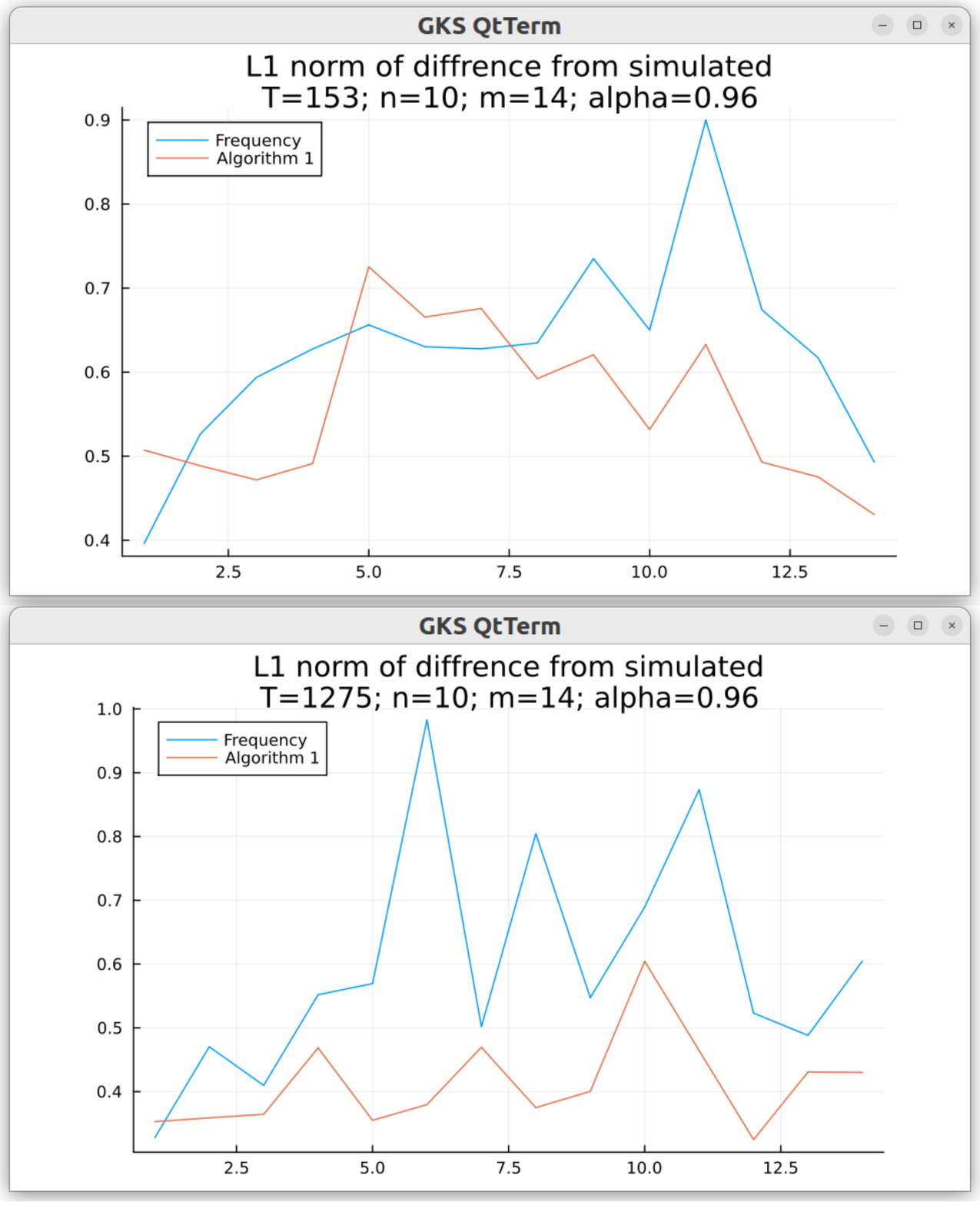}

	\caption{The top plot shows $l_1$-norm differences $\Delta_k(X)$ and $\Delta_k(R)$ for 14 daily random time intervals of lengths $14, 11, 12, 14, 11, 12, 9, 12, 7, 9, 11, 10, 9, 12$. Average $l_1$ norm difference is about $0.55$ for Algorithm 1 versus about $0.62$ for frequency counter. \\\hspace{\textwidth}\\\hspace{\textwidth}	
	The bottom plot shows $l_1$-norm differences $\Delta_k(X)$ and $\Delta_k(R)$ for 14 daily random time intervals of lengths $114, 89, 104, 94, 114, 115, 88, 99, 76, 37, 78, 72, 75, 120$. Average $l_1$ norm difference is about $0.40$ for Algorithm 1 versus about $0.59$ for frequency counter.  }
	\label{fig:screenshot}
\end{figure}

\noindent In the way just indicated, the simulation generates two categorical distribution sequences
\begin{enumerate}
	\item[] $X_t$ for the algorithm described in Section 2
	\item[]   a click frequency distribution $\hat{R_t} $ defined by averaging click counters 
	\begin{equation}
		\hat{R_t} = \left( \frac{ 1}{\sum_{i=1}^n  R_{t,i}} \right) R_t  \nonumber
	\end{equation} 
\end{enumerate}

\noindent The Algorithm 1 can be evaluated by comparing categorical distribution sequences $ X_t, \; \hat{R_t} $  to the "actual" incoming  click-distribution  sequences $C_k, \; k = 1,2, \cdots $. For example, one can compare average $l_1$-norm differences 
\begin{eqnarray}
\Delta_k(X)=  \frac{1}{|T_k|}\sum_{t \in T_k}  \parallel X_t - C_k \parallel_1    \nonumber \\ 
\Delta_k(R)=  \frac{1}{|T_k|} \sum_{t \in T_k}  \parallel\hat{R_t} - C_k \parallel_1 \nonumber 
\end{eqnarray}
The sequences $\Delta_k(X) ,
\Delta_k(R) $ for typical  simulation runs  are plotted (cf. \cite{Plots}) in Figures 1 above.

\end{document}